\documentclass{amsart}[12pt]
\usepackage{amsmath}

\newtheorem{theorem}{Theorem}[section]

\newtheorem{lemma}[theorem]{Lemma}

\theoremstyle{definition}

\theoremstyle{remark}

\newtheorem*{remarks*}{Remarks}

\numberwithin{equation}{section}

\renewcommand{\AA}{{\mathbb A}}

\newcommand{\ZZ}{{\mathbb Z}}

\newcommand{\RR}{{\mathbb R}}

\DeclareMathOperator{\const}{Const.}

\newcommand\Cinf{\mathcal{C}^\infty}

\newcommand{\st}[1]{\ensuremath{^{\scriptstyle \textrm{#1}}}}

\renewcommand{\span}{\mbox{\textrm span}}

\newcommand{\supp}{\mbox{\textrm supp}}

\begin{document}

\openup3pt

\title{Riemann sums over polytopes}
\author{V. Guillemin}\address{Department of Mathematics\\
MIT \\Cambridge, MA 02139}
\author{S. Sternberg}
\address{Department of Mathematics\\
Harvard University\\Cambridge, MA  02140}
\email{shlomo@math.harvard.edu}
\maketitle


\section{Introduction}
\label{sec:1}
Given a $\Cinf$ function, $f$, on the interval $[0,1]$ let $R_N (f)$
 be the Riemann sum

\begin{equation}
   \label{eq:1.1}
   \frac{1}{N}\sum^N_{i=1} f(t_i) \, , \quad
        \frac{i}{N}\leq t_i < \frac{i+1}{N}\, .
 \end{equation}

In freshman calculus one learns that

\begin{equation}
  \label{eq:1.2}
  R_N (f) = \int^1_0 f(x) \, dx  + O \left( \frac{1}{N}\right)\, .
\end{equation}
What is not perhaps as well known is that if one chooses the
$t_i$'s judiciously, i.e.,~lets $t_i = \frac{i}{N}$ the $O \left(
  \frac{1}{N}\right)$ in (\ref{eq:1.2}) can be replaced by a much
better error term, an asymptotic series:
\begin{equation}
  \label{eq:1.3}
 \frac{1}{2N} \left(f(1) - f(0)\right) + \sum^{\infty}_{k=1} (-1)^{k-1}
  \frac{B_k}{(2k)!} \left(f^{\langle 2k-1 \rangle} (1)-
 f^{{\langle 2k-1 \rangle}}(0) \right) N^{-2k}
\end{equation}
in which the $B_k$'s are the Bernoulli numbers.  In particular if
$f$ is periodic of period~1 the $O \left( \frac{1}{N}\right)$ in
(\ref{eq:1.2}) is actually an $O (N^{-\infty})$.  (For an
expository account of this ``Euler--Maclaurin formula for Riemann
sums'' see \cite{GS}.)

In this article we will prove an $n$-dimensional version of this
result in which the interval $[0,1]$ gets replaced by a convex
polytope.  We will give a precise formulation of our result in
\S\ref{sec:4}; however, roughly speaking, it asserts that if
$\Delta$ is a simple convex polytope whose vertices lie on the
lattice, $\ZZ^n$, and if $f$ is in $\Cinf (\Delta )$ the
difference
\begin{equation}
  \label{eq:1.4}
  \int_{\Delta} f(x) \, dx -\frac1{N^n}\sum_{k \in N \Delta \cap \ZZ^n}
     f \left( \frac{k}{N}\right)
\end{equation}
can be expanded in an asymptotic series in $N^{-1}$ in which the
coefficients are explicitly computable by recipes resembling
(\ref{eq:1.3}).  Our formula bears a formal resemblance to the
generalized Euler--Maclaurin formulas of \cite{KP}, \cite{KK},
\cite{CS}, \cite{Gu}, \cite{BV} et al., however in these so-called
``exact'' Euler--Maclaurin formulas the functions involved are
polynomials, not as in the case here, arbitrary $\Cinf$
functions.  Somewhat closer in spirit to our result is the
Euler--Maclaurin formula with remainder of \cite{KSW} and the
Ehrhard theorem for symbols of \cite{GSW}.  (Our result also
yields an Ehrhard theorem for symbols, and its relation to the
theorem in \cite{GSW} will be discussed in \S~\ref{sec:4}.)

A word about the organization of this paper.  In \S~\ref{sec:2}
 we will review the proof of the Riemann sum version of
 Euler--Maclaurin, for the interval, $(-\infty , 0]$ and in
 \S\ref{sec:3} show how to extend this result to regions in
 $\RR^n$ which are defined by systems of $k$ linearly-independent
 inequalities

 \begin{equation}
   \label{eq:1.5}
 \langle u_i , x \rangle \leq c_i \, , \quad
   u_i \in \ZZ^n \, , \quad
   c_i \in \ZZ \, .
 \end{equation}
(We will call such regions \emph{$k-wedges$}.)

In \S\ref{sec:4} we will derive from this result a
Euler--Maclaurin formula for Riemann sums over polytopes and in
\S\ref{sec:5} show that our result has an equivalent formulation
as an Ehrhard theorem for symbols.

We would like to thank Dan~Stroock and Hans Duistermaat for
helpful discussions concerning the material in
Section~\ref{sec:2}.

\section{Euler--Maclaurin for the interval $(-\infty , 0]$}

\label{sec:2}

Let $\tau (s)$ be the Todd function

\begin{equation}
  \label{eq:2.1}
 \frac{s}{1-e^{-s}} = 1 + \frac{s}{2} + \sum (-1)^{n-1}
B_n \frac{s^{2n}}{(2n)!}\, .
\end{equation}

In this section we will show that for Schwartz functions, $f \in
S(\RR)$ the difference

\begin{equation}
\label{eq:2.2}
 \frac{1}{N} \sum^{\infty}_{k=0} f \left( -\frac{k}{N}\right)
 - \int^0_{-\infty} f(x) \, dx
\end{equation}
has an asymptotic expansion:

\begin{equation}
  \label{eq:2.3}
    \frac{f(0)}{2N} + \sum^{\infty}_{n=1} (-1)^{n-1}
  \frac{B_n}{(2n)!} f^{(2n-1)}(0) N^{-2n} \, .
\end{equation}

In view of (\ref{eq:2.1}) this formula can be written more
succinctly in the form
\begin{equation}
 \label{eq:2.4}
 \frac{1}{N} \sum^{\infty}_{k=0} f \left( -\frac{k}{N}\right)
  \sim \left( \tau \left(\frac{1}{N}\, \frac{\partial}{\partial h}\right)
 \int^h_{-\infty} f(x) \, dx \right) (h=0)
\end{equation}
and it is this version of it which we will prove.

We first of
all observe that if $f(x) = e^{\lambda x}$, $\lambda >0$, then

\begin{displaymath}
  \int^{h}_{-\infty} f(x) \, dx = \frac{1}{\lambda} e^{\lambda h}\,.
\end{displaymath}
So for $N>2\pi \lambda$ we may apply the infinite order constant coefficient operator
$ \tau \left( \frac{1}{N} \, \frac{\partial}{\partial h}\right)$ to this expression:
\begin{eqnarray*}
 \tau \left( \frac{1}{N} \, \frac{\partial}{\partial h}\right)
  \int^h_{-\infty} f(x) \, dx 
 &=& \tau  \left( \frac{1}{N} \, \frac{\partial}{\partial  h}\right)
  \, \frac{e^{\lambda h}}{\lambda}\\[1ex]
 &=& \tau \left( \frac{\lambda}{N} \right) \, 
 \frac{e^{\lambda h}}{\lambda}\\[1ex]
&=& \frac{1}{N} \, \frac{\lambda}{1-e^{-\lambda /N}} \, 
\frac{e^{\lambda h}}{\lambda} \\[1ex]
&=& \frac{1}{N} \left( \sum^{\infty}_{k=0}
e^{-\frac{k}{N}\lambda} \right) e^{\lambda h},
\end{eqnarray*}
all series being convergent. 
We conclude that
\begin{equation}
  \label{eq:2.5}
 \frac{1}{N} \sum^{\infty}_{k=0} e^{-\frac{k}{N}\lambda}
= \left( \tau \left( \frac{1}{N}\, \frac{\partial}{\partial h}
 \right) \int^{h}_{-\infty} e^{\lambda x}\, dx \right)
(h=0)\, .
\end{equation}
More generally differentiating this identity $n$~times with
respect to $\lambda$ we obtain
\begin{equation}
 \label{eq:2.6}
\frac{1}{N} \sum^{\infty}_{k=0} \left(-\frac{k}{N}\right)^{n}
e^{-\frac{k}{N}\lambda} = \left( \tau \left( \frac{1}{N}
 \, \frac{\partial}{\partial h} \right) \int^h_{-\infty}
 x^n e^{\lambda x}\, dx \right) (h=0)
\end{equation}
verifying (\ref{eq:2.4}) for the function $x^n e^{\lambda x}$ and
hence for the functions of the form $p(x) e^{\lambda x}$ where
$p$ is a polynomial.  Now let $f$ be a Schwartz function and $p$
a polynomial having the property that $f(x) - p(x) e^{\lambda x}$
vanishes to order $n+2$ at $x=0$.  Let

\begin{equation}
  \label{eq:2.7}
 g(x) = 
\begin{cases}
0 \, , & x \geq 0 \\
f(x) -p(x) e^{\lambda k} \, , & x<0 \, .
\end{cases}
\end{equation}

Then

\begin{equation}
  \label{eq:2.8}
   \| g^{(i)} (x) \|_1 \leq \infty \hbox{  for  } i \leq n+2
\end{equation}
and by the Poisson summation formula
\begin{equation}
\label{eq:2.9}
 \sum_{-\infty <k<\infty} g \left(- \frac{k}{N}\right)
 = N \sum_{-\infty <k<\infty}\hat{g} (Nk) \, .
\end{equation}
However, by (\ref{eq:2.8})
\begin{equation}
  \label{eq:2.10}
  | \hat{g} (Nk) | \leq \const N^{-n}k^{-2}
\end{equation}
for $k \neq 0$, and
\begin{equation}
  \label{eq:2.11}
  \hat{g} (0) = \int^0_{-\infty} g (x) \, dx \, .
\end{equation}

Hence

\begin{equation}
  \label{eq:2.12}
  \frac{1}{N} \sum^{\infty}_{k=0} g \left( -\frac{k}{N}\right)
    = \int^0_{-\infty} g(x) \, dx +O(N^{-n}) \, .
\end{equation}

This shows that (\ref{eq:2.4}) is true for $g$ modulo $O(N^{-n})$
and hence is true for $f$ modulo $O (N^{-\infty})$.

\hfill Q.E.D

In \S\ref{sec:3} we will also need a version of the theorem above
for ``twisted'' Riemann sums.  Let $\omega\neq 1$ be a $q$\st{th}
 root of unity and let

  \begin{displaymath}
    \tau_{\omega} (s) = \frac{s}{1-\omega e^{-s}} 
       = \frac{s}{1-\omega} + \sum_{i>1} b^{\omega}_i s^i \, .
  \end{displaymath}

For $f \in S(\RR)$ we will show that the twisted Riemann sum

\begin{equation}
  \label{eq:2.13}
  \frac{1}{N} \sum^{\infty}_{k=0} \omega^k f \left( -\frac{k}{N}\right)
\end{equation}
is asymptotic to the series

\begin{equation}
  \label{eq:2.14}
  \frac{1}{1-\omega} \frac{f(0)}{N} + \sum_{i>1} b^{\omega}_i
      f^{(i)} (0) N^{-i}\, .
\end{equation}

As above we can rewrite this in the more succinct form

\begin{equation}
  \label{eq:2.15}
  \frac{1}{N} \sum^{\infty}_{k=0} \omega^k f 
      \left(-\frac{k}{N}\right) \sim
        \left( \tau_{\omega}\left( \frac{1}{N}\, 
            \frac{\partial}{\partial h}\right) \int^h_{-\infty}
            f(x) \, dx \right) (h=0)
\end{equation}
and we will prove this by essentially the same proof as before:
If $f = e^{\lambda x}$ the expression in parentheses is

\begin{eqnarray}
  \label{eq:2.16}
      \tau_{\omega} \left( \frac{\lambda}{N}\right) 
         \frac{e^{\lambda h}}{\lambda }
     &=& \frac{1}{N} 
          \left( \frac{\lambda}{1-\omega e^{-\lambda  /N}}\right)
          \frac{e^{\lambda h}}{\lambda}\\[1ex]\notag
     &=& \frac{1}{N} \left( \sum^{-\infty}_{k=0} \omega^k
           e^{-k \lambda /N} \right) e^{\lambda h } \, ,\notag
\end{eqnarray}
and by setting $h=0$ we see that (\ref{eq:2.15}) is valid for
$f=e^{\lambda x}$; and by differentiating both sides of
(\ref{eq:2.16}) by $\left( \frac{d}{d \lambda}\right)^n$ that
it's valid for $x^n e^{\lambda x}$ and hence for $p(x) e^{\lambda
x}$ where $p (x)$ is a polynomial.  Thus, as above, we're reduced
to showing that for the function $g$ defined by (\ref{eq:2.7}):

\begin{equation}
  \label{eq:2.17}
  \frac{1}{N}\sum_{-\infty <k<\infty} \omega^k g
  \left( \frac{k}{N} \right) =O (N^{-n})\, .
\end{equation}
For $r=0,1,\ldots,q-1$, let $g_r (x) = g(qx+ \frac{r}{N})$.  Then 

\begin{equation}
  \label{eq:2.18}
  \frac{1}{N} \sum_{-\infty < k < \infty} \omega^k g
  \left( \frac{k}{N}\right) = \frac{1}{N} \sum^{q-1}_{r=0}
     \omega^r \left( \sum_{-\infty < k < \infty} g_r
       \left( \frac{k}{N}\right) \right) \, .
\end{equation}

%




%

Since

\begin{displaymath}
  \hat{g}_r (Nk) = \frac{1}{q} \, e^{i\frac{rk}{q}} \hat{g} 
     \left(\frac{Nk}{q}\right)
\end{displaymath}
the Poisson summation formula yields, as before, the estimate

\begin{equation}
  \label{eq:2.19}
  \sum^{q-1}_{r=0} \omega^r \int^{\infty}_{-\infty}
     g_r (x) \, dx + O (N^{-n})\, ,
\end{equation}
for the right hand side of (\ref{eq:2.18}).  However,

\begin{displaymath}
  \int^{\infty}_{-\infty} g_r (x) \, dx 
      = \int^{\infty}_{-\infty} g_0 (x) \, dx
\end{displaymath}
and $\sum^{q-1}_{r=0} \omega^r =0$ so the first summand in
(\ref{eq:2.18}) is zero.
\hfill Q.E.D.

We will conclude this discussion of one dimensional
Euler--Maclaurin formulas by describing analogues of
(\ref{eq:2.4}) and (\ref{eq:2.15}) in which the sum over $-\infty
< k<0$ gets replaced by a sum over $-\infty < k < \infty$.  For
simplicity assume that $f \in \Cinf_0 (\RR)$.  We claim:

\begin{equation}
  \label{eq:2.20}
  \frac{1}{N} \sum^{\infty}_{k=-\infty} f 
      \left( \frac{k}{N}\right) = \int^{\infty}_{-\infty}
        f(x) \, dx + O (N^{-\infty})
\end{equation}
and, for $\omega$ a $q$\st{th} root of unity, $\omega\neq 1$,
\begin{equation}
  \label{eq:2.21}
  \sum^{\infty}_{k=-\infty} \omega^k f \left( \frac{k}{N}\right)
      = O (N^{-\infty})\, .
\end{equation}
To prove (\ref{eq:2.20}) we first observe that for $c$ a large
positive integer, the left and right hand sides of
(\ref{eq:2.20}) are unchanged if one substitutes the function, $f
(x+c)$, for $f$, so without loss of generality we can assume that
$f$ is supported on the interval, $x<0$, in which case
(\ref{eq:2.3}) is of order $O (N^{-\infty})$ and (\ref{eq:2.20})
is a consequence of (\ref{eq:2.4}).  Similarly if we replace
$f(x)$ by $f(x+cq)$, with $c$ a large positive integer, the left
and right hand sides of (\ref{eq:2.21}) are unchanged; so we can
assume that $f$ is supported on the interval $x<0$, and
(\ref{eq:2.21}) is a consequence of (\ref{eq:2.15}).

\section{Euler--Maclaurin for wedges}
\label{sec:3}

Let $\ZZ^n$ be the integer lattice in $\RR^n$, $(\ZZ^n)^*$ its
dual lattice in $(\RR^n)^*$ and $\langle u,x \rangle$ the usual
paring of vectors, $x \in \RR^n$, and $u \in (\RR^n)^*$.  Given
$m$ linearly independent vectors, $u_i \in (\RR^n)^*$ we will
call the subset of $\RR^n$ defined by the inequalities
\begin{equation}
  \label{eq:3.1}
        \langle u_i,x \rangle \leq c_i \qquad i=1,\ldots ,m                  
\end{equation}
an \emph{integer $m$-wedge} if the $c_i$'s are integers and the
$u_i$'s primitive lattice vectors in $(\ZZ^n)^*$.  Let $W$ be the
set (\ref{eq:3.1}) and $U$ the subspace of $(\RR^n)^*$ spanned by
the $u_i$'s.  We will call $W$ a \emph{regular} integer $m$-wedge
  if $u_1,\ldots ,u_m$ is a lattice basis of the lattice $U \cap
(\ZZ^n)^*$ i.e.,~if
\begin{equation}
  \label{eq:3.2}
  U \cap (\ZZ^n)^* = \span_{\ZZ} \{ u_1,\ldots ,u_m \}\, .
\end{equation}
We will need below the following criterion for regularity.

\begin{lemma}
  \label{lem:3.1}
If (\ref{eq:3.2}) holds, $u_i,\ldots ,u_m$ can be extended to a
lattice basis, $u_1 ,\ldots ,u_n$ of $(\ZZ^n)^*$.
\end{lemma}
\begin{proof}
Let $u_{m+1},\ldots ,u_n$ be vectors in $(\ZZ^n)^*$ whose
projections onto the quotient of $(\ZZ^n)^*$by $U \cap (\ZZ)^*$
are a lattice basis of this quotient.
\end{proof}

For an integer $m$-wedge satisfying (\ref{eq:3.2}) the
$n$-dimensional generalization of Euler--Maclaurin is relatively
straightforward.
\begin{theorem}
  \label{th:3.2}
Let $W_h$ be the subset of $\RR^n$ defined by the inequalities
\begin{equation}
  \label{eq:3.3}
  \langle u_i,x \rangle \leq c_i + h_i \, , \quad
  i=1,\ldots ,m \, .
\end{equation}
\end{theorem}

Then, for $f \in \Cinf_0 (\RR^n)$,
\begin{equation}
  \label{eq:3.4}
  \frac{1}{N^n} \sum_{k \in \ZZ^n \cap NW} f 
     \left(\frac{k}{N}\right) \sim \left( \tau \left(
         \frac{1}{N}\, \frac{\partial}{\partial h}\right)
         \int_{W_h} f(x) \, dx \right) (h=0)
\end{equation}
where $\tau (s_1 ,\ldots ,s_m) = \prod^m_{i=1} \tau (s_i)$.
\begin{proof}
By Lemma~\ref{lem:3.1} we can incorporate $u_1,\ldots ,u_m$ in a
lattice basis $u_1 ,\ldots ,u_n$ of $(\ZZ^n)^*$.  Let
$\alpha_1,\ldots ,\alpha_n$ be the dual basis of $\ZZ^n$ and let
$v = \sum^m_{i=1}c_i\alpha_i$.  Then via the map
\begin{equation}
  \label{eq:3.5}
  x \in \RR^n \to \sum x_i \alpha_i +v
\end{equation}
one is reduced to proving the theorem for the standard $m$-wedge:
 $x_1 \leq 0 , \ldots , x_m \leq 0$, i.e.,~showing that the sum
 \begin{equation}
   \label{eq:3.6}
   \frac{1}{N^n} \sum f \left( \frac{k_1}{N},\cdots , 
        \frac{k_n}{N}\right)
 \end{equation}
summed over all $(k_1, \ldots ,k_n) \in \ZZ^n$, with $k_i \leq 0$
for $i \leq m$, is equal to the expression
\begin{equation}
  \label{eq:3.7}
  \tau \left( \frac{1}{N}\, \frac{\partial}{\partial h}\right)
    \int^{h_1}_{-\infty} \cdots \int^{h_m}_{-\infty}
      dx_1 \ldots dx_m \int^{\infty}_{-\infty} \cdots
      \int^{\infty}_{-\infty} f(x) \, 
         dx_{m+1} \cdots  dx_{n}\, ,
\end{equation}
evaluated at $h=0$, modulo $O (N^{-\infty})$.  Moreover, since
the subalgebra of $\Cinf_0 (\RR^n)$ generated by the products 
\begin{displaymath}
  f(x) = f_1(x_1) \ldots f_n (x_n) \, , \quad f_i \in \Cinf_0 (\RR)\, ,
\end{displaymath}
is dense in $\Cinf_0 (\RR^n)$ it suffices to prove the theorem
for functions of this form, and hence it suffices to prove the
theorem for $n=1$ and $m=0$ or $1$.  However, these two cases
were dealt with in \S\ref{sec:2}.  (See~(\ref{eq:2.4}) and
(\ref{eq:2.20}).)
\end{proof}

We will next describe how (\ref{eq:3.4}) has to be modified if
the condition (\ref{eq:3.2}) isn't satisfied.  As  above let
$u_{m+1},\ldots , u_n$ be vectors in   $(\ZZ^n)^*$ whose
projections onto the quotient of $(\ZZ^n)^*$ by $U \cap (\ZZ^n)^*$
are a lattice basis of this quotient lattice.  The vectors,
$u_1,\ldots ,u_n$ are now no longer a lattice basis of
$(\ZZ^n)^*$ but they span a sublattice

\begin{equation}
\label{eq:3.8}
 \AA^* = \span_{\ZZ} \{ u_1,\ldots, u_n \}                                
\end{equation}
of $(\ZZ^n)^*$ of rank~$n$, so the quotient
\begin{equation}
  \label{eq:3.9}
  \Gamma = (\ZZ^n)^* /\AA^*
\end{equation}
is a finite group.  Let $\alpha_1,\ldots ,\alpha_n$ be the basis
vectors of $\RR^n$ dual to $u_1,\ldots ,u_n$.  Since $\AA^*$ is a
sublattice of $(\ZZ^n)^*$ the dual lattice,
\begin{equation}
  \label{eq:3.10}
  \AA = \span_{\ZZ} \{ \alpha_1 ,\ldots , \alpha_n \}\, ,
\end{equation}
contains $\ZZ^n$ as a sublattice.  Moreover, each element, $x \in
\AA$, defines a character of the group, $\Gamma$, via the pairing
\begin{equation}
  \label{eq:3.11}
  \gamma \in \Gamma \to e^{2\pi i \langle \gamma , x \rangle}
\end{equation}
and this character is trivial if and only if $x$ is in $\ZZ^n$.
By a theorem of Frobenius the average value of a character of a
finite group is zero if the character is non-trivial and is one
if it is trivial, so we have
\begin{equation}
  \label{eq:3.12}
  \frac{1}{|\Gamma |} \sum e^{2\pi i \langle \gamma ,x \rangle}=
\begin{cases}
1 \hbox{  if  } x \in \ZZ^n\\
0 \hbox{  if  } x \notin \ZZ^n \, .
\end{cases}
\end{equation}
For each $\gamma \in \Gamma$ let
\begin{equation}
  \label{eq:3.13}
  \tau_{\gamma} (s_1,\ldots ,s_m) 
      = \tau_{\omega_1}(s_1)  \cdots \tau_{\omega_m} (s_m)
\end{equation}
where $\omega_k = e^{2\pi i \langle \gamma ,\alpha_{k} \rangle}$.  We will
generalize Theorem~\ref{th:3.2} by showing that for integer
$m$-wedges which don't satisfy condition~(\ref{eq:3.2}) one has
\begin{theorem}
  \label{th:3.3}
For $f \in \Cinf_0 (\RR^n)$
\begin{equation}
  \label{eq:3.14}
  \frac{1}{N^n} \sum_{k \in N W \cap \ZZ^n}
  f \left( \frac{k}{N}\right) = \left( \sum_{\gamma \in \Gamma}
    \tau_{\gamma} \left( \frac{1}{N}\, \frac{\partial}{\partial h} 
    \right)
    \int_{W_h} f(x) \, dx \right) (h=0) \mod O (N^{-\infty})\, .
\end{equation}
\end{theorem}
\begin{proof}
By (\ref{eq:3.11}) the sum on the left coincides with the sum
\begin{equation}
  \label{eq:3.15}
  \frac{1}{|\Gamma |} \sum_{\gamma \in \Gamma} \frac{1}{N^n}
  \sum_{x \in \AA \cap NW} e^{2\pi i \langle \gamma ,x \rangle}
  f \left( \frac{x}{N}\right)
\end{equation}
so it suffices to show that the $\gamma$-\st{th} summand in
(\ref{eq:3.14}) is equal to the $\gamma$-\st{th} summand in
(\ref{eq:3.15}).  Via the map~(\ref{eq:3.5}) the $\gamma$-\st{th}
summand in (\ref{eq:3.15}) becomes
\begin{equation}
  \label{eq:3.16}
  \frac{1}{N^n |\Gamma |} \sum_{k_1 \leq 0 ,\ldots ,k_m \leq 0}
     \omega_1^{k_1} \ldots \omega^{k_m}_m \left(
\sum_{k_{m+1, \ldots ,k_n}} g \left( \frac{k}{N}\right)\right)
\end{equation}
where $g (x_1,\ldots ,x_n) = f (v + x_1\alpha_1 + \cdots +
x_n\alpha_n)$, and the  $\gamma$-\st{th} summand in
(\ref{eq:3.14}) becomes 

\begin{equation}
  \label{eq:3.17}
  \frac{1}{|\Gamma |}\,  \tau_{\gamma} 
     \left( \frac{1}{N}\, \frac{\partial}{\partial h}\right)
       \int^{h_1}_{-\infty} \ldots \int^{h_m}_{-\infty}
         dx_1 \ldots dx_m \int^{\infty}_{-\infty}\ldots
           \int^{\infty}_{-\infty} g(x) dx_{m+1}\ldots dx_k
\end{equation}
evaluated at $h=0$.  (The reason for the factor, $1/|\Gamma |$,
is that this is the Jacobian determinant of the mapping
(\ref{eq:3.5}).)  To prove that (\ref{eq:3.16}) and
(\ref{eq:3.17}) are equal mod $O (N^{-\infty})$ it suffices as
above to prove this for functions of the form $g=g_1(x) \ldots g_n(x_n)$
 with $g_i \in \Cinf_0 (\RR)$ and hence to show, for $i \leq m$
 \begin{equation}
   \label{eq:3.18}
   \frac{1}{N} \sum^{-\infty}_{k_i=0}\omega^{k_i} g_i
   \left( \frac{k_i}{N} \right) \sim \left( \tau_{\omega_i}
     \left( \frac{1}{N}\, \frac{\partial}{\partial h_i} \right)
     \int^{h_i}_{-\infty} g_i (x_i)\, dx_i \right)(h_i=0)
 \end{equation}
and, for $i >m$
\begin{equation}
  \label{eq:3.19}
  \frac{1}{N} \sum^{\infty}_{-\infty} g_i 
    \left( \frac{k_i}{N}\right) = \int^{\infty}_{-\infty} g_i
      (x_i) \, dx_i + O(N^{-\infty})\, ,
\end{equation}
and these follow from the identities (\ref{eq:2.15}) and (\ref{eq:2.20}).
\end{proof}

\section{Riemann sums over polytopes}
\label{sec:4}

Let $\Delta \subseteq \RR^n$ be an $n$-dimensional polytope whose
verticies lie on the lattice $\ZZ^n$.  $\Delta$ is said to be a
simple polytope if each codimension~$k$ face is the intersection
of exactly $k$ facets.  (It suffices to assume that the vertices
of $\Delta$, i.e.,~the codimension $n$ faces, have this property
or, alternatively, that there are exactly $n$ edges of $\Delta$
meeting at each vertex.)  If the number of facets is $d$ then
$\Delta$ can be defined by a set of $d$ inequalities
\begin{equation}
  \label{eq:4.1}
  \langle u_i,x \rangle \leq c_i
\end{equation}
where $c_i$ is an integer and $u_i \in (\ZZ^n)^*$ is a primitive
lattice vector which is perpendicular to the $i$\st{th} facet
and points `` outward'' from $\Delta$.  By the simplicity
assumption each codimension $k$ face of $\Delta$ is the
intersection of $k$ facets lying in the hyperplanes
\begin{equation}
  \label{eq:4.2}
   \langle u_i,x \rangle = c_i\, \quad i \in F
\end{equation}
where $F$ is a $k$ element subset of $\{ 1,\ldots ,d \}$.  Let
$W_F$ be the $k$-wedge
\begin{equation}
  \label{eq:4.3}
   \langle u_i,x \rangle \leq c_i\, \quad i \in F\, .
\end{equation}

We will say that $\Delta$ is regular if each of these $k$-wedges
is regular.  (As above it suffices to assume this for the zero
faces, i.e.,~the vertices of $\Delta$, or alternatively to assume
that for every vertex, $v$, the edges of $\Delta$ which
intersect at $v$ lie on rays
\begin{displaymath}
  v + t \alpha_i \, , \quad 0 \leq t <\infty
\end{displaymath}
where $\alpha_1 ,\ldots ,\alpha_n$ is a lattice basis of
$\ZZ^n$.)

For regular simple lattice polytopes one has the following
Euler--Maclaurin formula.
\begin{theorem}
  \label{th:4.1}
Let $\Delta_h$ be the polytope
\begin{equation}
  \label{eq:4.4}
   \langle u_i,x \rangle \leq c_i + h_i \quad i=1,\ldots ,d \, .
\end{equation}
Then, for $f \in \Cinf_0 (\RR^n)$
\begin{equation}
  \label{eq:4.5}
  \frac{1}{N^n} \sum_{k \in \ZZ^n \cap N\Delta} f 
     \left( \frac{k}{N}\right) \sim \left( \tau \left(
           \frac{1}{N} \, \frac{\partial}{\partial h}\right)
           \int_{\Delta_h} f(x) \, dx \right) (h=0)
\end{equation}
where $\tau (s_1,\ldots ,s_d) = \tau (s_1) \ldots \tau (s_d)$.
\end{theorem}
\begin{proof}
By a partition of unity argument we can assume that $\supp f$ is
contained in a small neighborhood of the set (\ref{eq:4.2})
and doesn't intersect the hyperplanes, $ \langle u_i,x \rangle =
c_i$, $i \notin F$.  Then for $i \notin F$
\begin{displaymath}
  \tau \left( \frac{1}{N} \, \frac{\partial}{\partial h_i}\right)
      \int_{\Delta_h} f \, dx = \int_{\Delta_h} f \, dx 
        + \frac{1}{2N} \,  \frac{\partial}{\partial h_i}
          \int_{\Delta_h} f(x) \, dx + \cdots \, .
\end{displaymath}
\end{proof}
However, by (\ref{eq:2.4}) all the terms on the right except the first are
integrals of derivatives of $f$ over the hyperplane  $ \langle u_i,x \rangle =
c_i +h_i$, and hence for $h_{i}$ small are zero.  Thus the left
hand side of (\ref{eq:4.5}) becomes
\begin{displaymath}
  \left( \prod_{i \in F} \tau \left( \frac{1}{N}\, 
      \frac{\partial}{\partial h_i} \right) \int_{(W_F)_h}
      f(x) \, dx \right) (h=0)
\end{displaymath}
and the theorem above reduces to Theorem~\ref{th:3.2}.
\hfill Q.E.D.

If $\Delta$ is simple but not regular, one gets a slightly more complicated
result.  To the codimension $k$-face of $\Delta$
 defined by (\ref{eq:4.2}) attach the subspace
 \begin{displaymath}
   U_F = \span_{\RR} \{ u_i \, , \quad i \in F \}
 \end{displaymath}
of $(\RR^n)^*$, the sublattice
\begin{displaymath}
  \ZZ_{\Gamma} = \span_{\ZZ} \{  u_i \, , \quad i \in F \}
\end{displaymath}
and the finite group
\begin{displaymath}
  \Gamma_F = U_F \cap (\ZZ^n)^* /\ZZ_F \, .
\end{displaymath}
This group coincides with the ``torsion group'' (\ref{eq:3.9}) of
the wedge $W_F$.  Moreover, if $E$ is a subset of $F$, $U_E$ is
contained in $U_F$ and $\ZZ_E$ in $\ZZ_F$, so $\Gamma_E$ is
contained in $\Gamma_F$.  Let $\Gamma^{\sharp}_F$ be the set of
points in $\Gamma_F$ which are \emph{not} contained in $\Gamma_E$
for some proper subset, $E$ of $F$.
\begin{theorem}
  \label{th:4.2}
For $f \in \Cinf_0 (\RR^n)$ the sum
\begin{equation}
  \label{eq:4.6}
\frac1{N^n}  \sum_{k \in \ZZ^n \cap N \Delta} f\left( \frac{k}{N}\right)
\end{equation}
is equal $\mod O (N^{-\infty})$ to 
\begin{equation}
  \label{eq:4.7}
  \left( \sum_F \sum_{\gamma \in \Gamma^{\sharp}_F} \tau_{\gamma}\left(
      \frac{1}{N} \, \frac{\partial}{\partial h} \right)
        \int_{\Delta_h} f(x) \, dx \right) (h=0)\, .
\end{equation}
\end{theorem}
\begin{proof}
As above it suffices to prove this for $\supp f$ contained in a
small neighborhood of the set (\ref{eq:4.2}), and not intersecting the
hyperplanes, $\langle u_i,x \rangle =c_i$, $i \notin F$.  Then as
above, the only contribution to the sum (\ref{eq:4.7}) is
\begin{displaymath}
  \left( \sum_{\gamma \in \Gamma_F} \tau_{\gamma} \left(
      \frac{1}{N} \, \frac{\partial}{\partial h} \right)
    \int_{(W_F)_h} f(x) \, dx \right) (h=0)
\end{displaymath}
and Theorem~\ref{th:4.2} reduces to Theorem~\ref{th:3.3}.
\end{proof}

\section{An Ehrhart theorem for symbols}
\label{sec:5}
A function, $f \in \Cinf (\RR^n)$ is a \emph{polyhomogeneous
  symbol}  of degree~$d$ if, for large values of $x$, it admits
an asymptotic expansion
\begin{equation}
  \label{eq:5.1}
  f(x) \sim \sum^{-\infty}_{j=d} f_j (x)
\end{equation}
whose summands are homogeneous functions $f_j \in \Cinf (\RR^n -
\{ 0 \})$ of degree~$j$.  Let $f$ be such a function and let
$\Delta $ be a simple lattice polytope in $\RR^n$ containing the
origin in its interior.  The Ehrhart function of the pair, $f$,
$\Delta$, is defined to be the function
\begin{displaymath}
  E (f,\Delta , N) = \sum_{k \in N \Delta \cap \ZZ^n} f(k) \, , 
  \quad N \in \ZZ_+ \, .
\end{displaymath}
In \cite{GSW} it was shown that  
$$  E (f,\Delta , N) -\int_{N\Delta}fdx$$
had an asymptotic
expansion
\begin{equation}
  \label{eq:5.2}
  \sum^{-\infty}_{j=n+d} c_j N^j + c
\end{equation}
for $N$ large.

  The main result of this section
is a variant of this result.  As above let $\Delta$ be a simple
lattice polytope in $\RR^n$ and let $C_{\Delta}$ be the polyhedral
cone consisting of all points, $(x_1,\ldots ,x_n, x_{n+1})$, in
$\RR^{n+1}$ with $x_{n+1} >0$ and $(x_1 ,\ldots ,x_n)/x_{n+1} \in
\Delta$.  Then, for $N \in \ZZ_+$, $N \Delta$ is just the slice of
$C_{\Delta}$ by the hyperplane, $x_{n+1}=N$.  We will prove that
if $f \in \Cinf (\RR^{n+1})$ is a homogeneous symbol of
degree~$d$ the sum
\begin{equation}
  \label{eq:5.3}
  \sum_{k \in N\Delta \cap \ZZ^n} f(k)
\end{equation}
admits an asymptotic expression of the form~(\ref{eq:5.2}).
\begin{remarks*}

 \begin{enumerate}
  \item  
This result, albeit very close in spirit to the
    theorem in \cite{GSW} cited above, doesn't, as far as we can
    see, seem to be a trivial consequence of it.
\item   
This result has a number of applications to spectral
theory on toric varieties which we'll explore in future publications.
\item As a corollary of this result one gets another variant of the Ehrhart theorem: Let
$$\Delta^\sharp :=\{(x_1,\dots, x_{n+1})\in C_\Delta,\ \ x_{n+1}\leq 1\}.$$
This $(n+1)$-dimensional polytope is {\em not} in general simple. However a version of 
theorem described at the beginning of this section is still true, namely
$$E(f,\Delta^\sharp,N)\sim \sum_{i=d+n+1}^{-\infty} c_i^\sharp + c^\sharp \log N.$$
as one can see by summing the differences
$$E(f,\Delta^\sharp,N)-E(f,\Delta^\sharp,N-1)$$
and noting 
that each difference is exactly (\ref{eq:5.3}). By combining this result with the Danilov ``desingularization trick" \cite{Da} one can extend the Ehrhart theorem to a much larger class of convex lattice polytopes.
We will discuss the details elsewhere.
  \end{enumerate}
\end{remarks*}
\begin{proof}
As above let
\begin{displaymath}
  f \sim \sum^{-\infty}_{i=d} f_i
\end{displaymath}
where $f_i (x_1,\ldots ,x_{n+1})$ is a homogeneous function of
degree~$i$.  Then on the cone, $C_{\Delta}$:
\begin{displaymath}
 f_i (x_1 ,\ldots ,x_{n+1}) = x^i_{n+1} \, f_i \left(
     \frac{x_1}{x_{n+1}}\, , \cdots \, , 
     \frac{x_n}{x_{n+1}} \, ,  \, 1 \right)
\end{displaymath}
so if we set $\tilde{f}_i (x_1,\ldots ,x_n) = f_i (x_1,\ldots
,x_n,1)$ the sum~(\ref{eq:5.3}) is equal to the sum
\begin{equation}
  \label{eq:5.4}
  N^i \sum_{k \in N \Delta \cap \ZZ^n}
     \tilde{f}_i \left( \frac{k}{N}\right)\, .
\end{equation}
which is $N^{i+n}$ times the Riemann sum
\begin{equation}
  \label{eq:5.5}
  \frac{1}{N^n}\sum_{k \in N \Delta \cap \ZZ^n}
     \tilde{f}_i \left( \frac{k}{N}\right)\, .
\end{equation}
Thus, by Theorem~\ref{th:4.2}, each of these summands admits an
asymptotic expansion:
\begin{displaymath}
  \sum^{-\infty}_{k=n+i} c_{i,k} N^k
\end{displaymath}
and hence so does the sum (\ref{eq:5.3}).
\end{proof}

\end{document}